# Fractional constant elasticity of variance model*

## Ngai Hang Chan[1] and Chi Tim Ng

*The Chinese University of Hong Kong*

**Abstract:** This paper develops a European option pricing formula for fractional market models. Although there exist option pricing results for a fractional Black-Scholes model, they are established without accounting for stochastic volatility. In this paper, a fractional version of the Constant Elasticity of Variance (CEV) model is developed. European option pricing formula similar to that of the classical CEV model is obtained and a volatility skew pattern is revealed.

## 1. Introduction

Black [3, 4] developed the so-called Constant Elasticity of Variance model for stock price processes that exhibit stochastic volatility. The CEV model is expressed in terms of a stochastic diffusion process with respect to a standard Brownian motion

$$(1.1) \qquad dX_t = \mu X_t \, dt + \sigma X_t^{\beta/2} \, dB_t,$$

where $0 \leq \beta \leq 2$ is a constant. If $\beta = 2$, such a model degenerates to a geometric Brownian motion. This model is characterized by the dependence of the volatility rate, i.e., $\sigma X^{\beta/2}$ on the stock price. When the stock price increases, the instantaneous volatility rate decreases. This seems reasonable because the higher the stock price, the higher the equity market value, and thus the lower the proportion of liability, which results in a decrease in the risk of ruin. The volatility rate or the risk measure is thus decreased. Making use of methods proposed in an earlier literature [6], Cox [5] studied the CEV models and gave an option pricing formula that involves a noncentral $\chi^2$ distribution function.

The classical CEV model (1.1) does not account for long-memory behavior, however. There are some evidences showing that the financial market exhibits long-memory structures (see [7, 23]). To encompass both long-memory and stochastic volatility, a possible model is to replace the Brownian motion in the stochastic diffusion equation by a fractional Brownian motion that exhibits a long-memory dependence structure (see [2, 13, 30]).

Though fractional Brownian motion can be used to model long-memory, as pointed out by Rogers in [29], the fractional Brownian motion is not a semi-martingale and the stochastic integral with respect to it is not well-defined in the classical Itô's sense. A theory different from the Itô's one should be used to handle the fractional situation. One approach is white noise calculus (see [18, 22, 28]), which was used in [14, 19] to construct stochastic integral with respect to the fractional Brownian motion. With the white noise approach, an extension to the Black-Scholes' stochastic differential equation is proposed to cope with long-memory phenomena (see [26]).

*Research supported in part by HKSAR RGC grants 4043/02P and 400305.
[1]The Chinese University of Hing Kong, e-mail: `nhchan@sta.cuhk.edu.hk`

In this paper, the white noise calculus approach is adopted to construct a fractional CEV model and to derive the option pricing formula for European call option. Basic concepts in white noise calculus are briefly introduced in Section 2. The fractional Itô's lemma, which is fundamental to option pricing theory, is presented in Section 3. Section 4 explains under what circumstances is the Itô's lemma applicable. In Section 5, the fractional option pricing theory is introduced and the concept of self-financing strategy, which is different from the traditional definition adopted by, e.g., Delbaen [10, 11], will further be discussed. Finally, the pricing formula for fractional CEV model is given in Sections 6.

## 2. White noise calculus and stochastic integration

In this section, the concept of stochastic integration with respect to fractional Brownian motion is introduced briefly. Important concepts are defined based on the white noise theory originated from [17], who considered the sample path of a Brownian motion as a functional. Throughout this paper, notations used in [1, 14, 18, 22] are adopted.

Let $S(R)$ be the Schwarz space. Take the dual $\Omega = S'(R)$, equipped with the weak star topology, as the underlying sample space, i.e., $\omega \in \Omega$ is a functional that maps a rapidly decreasing function $f(\cdot) \in S(R)$ to a real number. Also, let $B(\Omega)$ be the $\sigma$-algebra generated by the weak star topology. Then according to Bochner-Minlos Theorem (see Appendix A of [18]), there exists a unique probability measure $\mu$ on $B(\Omega)$, such that for any given $f \in S(R)$, the characteristic function of the random variable $\omega \to \omega(f)$ is given by

$$\int_\Omega e^{i\omega(f)} \, d\mu(\omega) = e^{-\frac{1}{2}||f||^2},$$

where

$$||f||^2 = \int_R f^2(t) \, dt.$$

Let $L^2$ be the space of real-valued functions with finite square norm $||\cdot||$, we have the triple $S(R) \subset L^2 \subset S'(R)$. For any $f \in L^2$, we can always choose a sequence of $f_n \in S(R)$ so that $f_n \to f$ in $L^2$, and $\omega(f)$ is defined as the $\lim_{n\to\infty} \omega(f_n)$ in $L^2(\mu)$.

Consider the indicator function

$$1_{(0,a)}(s) = \begin{cases} 1, \text{if } 0 \leq s < a, \\ -1, \text{if } a \leq s < 0, \\ 0, \text{otherwise}. \end{cases}$$

It can be verified that for any two real numbers $a$ and $b$, the random variables $\omega(1_{(0,a)}(\cdot))$ and $\omega(1_{(0,b)}(\cdot))$ are jointly normal, mean zeros, and with covariance $\min(a,b)$. Define $\tilde{B}(t)$ as $\omega(1_{(0,t)}(\cdot))$, we can always find a continuous version of $\tilde{B}(t)$, denoted by $B$, which is the standard Brownian motion. Roughly speaking, the probability space $(\Omega, B(\Omega), \mu)$ can intuitively be considered as a space consisting of all sample paths of a Brownian motion.

Following the approach of [14], we give the definition of fractional Brownian motion with Hurst parameter $\frac{1}{2} \leq H < 1$ in terms of white noise setting by using the fundamental operator $M_H$, defined on the space $L^2$, by

$$M_H f(x) = \begin{cases} f(x), & H = \frac{1}{2}, \\ c_H \int_R f(t)|t-x|^{H-\frac{3}{2}} dt, & H < \frac{1}{2} < 1, \end{cases}$$

where $c_H$ is a constant depending on the Hurst parameter $H$ via

$$c_H = (2\Gamma(H-\frac{1}{2})\cos(\frac{\pi}{2}(H-\frac{1}{2})))^{-1}[\Gamma(2H+1)\sin(\pi H)]^{\frac{1}{2}}.$$

Then, $\omega(M_H 1_{(0,a)}(\cdot))$ and $\omega(M_H 1_{(0,b)}(\cdot))$ are jointly normal with covariance

$$\frac{1}{2}\{|a|^{2H} + |b|^{2H} - |a-b|^{2H}\}.$$

Again, we can find a continuous version of $\omega(M_H 1_{(0,t)}(\cdot))$ denoted by $B^H(t)$, which is the fractional Brownian motion.

We have the following Wiener-Itô Chaos Decomposition Theorem for a square integrable random variable on $S'(R)$ (see Theorem 2.2.4 of [18]).

**Theorem 2.1.** *If $F \in L^2(\Omega, B(\Omega), \mu)$, then $F(\omega)$ has a unique representation*

$$F(\omega) = \sum_\alpha c_\alpha H_\alpha(\omega),$$

*where $\alpha$ is any finite integers sequence $(\alpha_1, \alpha_2, \ldots, \alpha_n)$, $c_\alpha$ are real coefficients and $H_\alpha(\omega) = h_{\alpha_1}(\omega(e_1))h_{\alpha_2}(\omega(e_2))\cdots h_{\alpha_n}(\omega(e_n))$ $h_n(x)$ are Hermite polynomials and $e_n$ is an orthonormal set in $S(R)$ which is defined as*

$$e_i(t) = (\sqrt{\pi}2^{i-1}(i-1)!)^{-1/2} h_{i-1}(t) e^{-t^2/2}.$$

*Furthermore, the $L^2$ norm of the functional $F(\omega)$ is given by*

$$\sum_\alpha \alpha! c_\alpha^2,$$

*with $\alpha! = \alpha_1! \alpha_2! \cdots \alpha_n!$.*

**Remark 2.1.** (see [18, 27]) The basis $\{H_\alpha(\omega) : \alpha\}$ is orthogonal with respect to the inner product $E(XY)$ in $L^2(\Omega, B(\Omega), \mu)$. The variance of $H_\alpha(\omega)$ is $\alpha!$. $H_0(\omega)$ is taken as the constant 1. For $\alpha \neq 0$, the expectation of $H_\alpha(\omega)$ is

$$E(H_\alpha(\omega) H_0(\omega)) = 0.$$

As a result, the term $c_0$ is the expectation of the functional $F(\omega)$.

Consider the functional $B_t^H = \omega(M_H 1_{[0,t]}(.))$, where $t$ is a given constant. Using the dual property (see [14]) of the $M_H$ operator: for all rapidly decreasing functions $f$ and $g$, we have

$$(f, M_H g) \equiv \int_R f(t) M_H g(t)\, dt = \int_R g(t) M_H f(t)\, dt \equiv (M_H f, g).$$

The function $M_H 1_{[0,t]}(.)$ can be rewritten in the Fourier expansion form:

$$M_H 1_{[0,t]}(s) = \sum_{i=0}^{\infty} (M_H 1_{[0,t]}(\cdot), e_i(\cdot)) e_i(s)$$

$$= \sum_{i=0}^{\infty} (1_{[0,t]}(\cdot), M_H e_i(\cdot)) e_i(s)$$

$$= \sum_{i=0}^{\infty} \left\{ \int_0^t M_H e_i(u) du \right\} e_i(s).$$

Since $\omega$ is linear, the functional can be written as

$$\omega(M_H 1_{[0,t]}(\cdot)) = \sum_{i=0}^{\infty} (1_{[0,t]}(\cdot), e_i(\cdot)) \omega(e_i)$$

$$= \sum_{i=0}^{\infty} \left\{ \int_0^t M_H e_i(u) du \right\} \omega(e_i).$$

$$= \sum_{\alpha} c_\alpha H_\alpha(\omega).$$

In this example, when $\alpha = \epsilon(i) = \{0..., 1, 0, ...\}$, i.e., one at position $i$, $c_\alpha = \int_0^t M_H e_i(u) du$, and $c_\alpha = 0$ otherwise. It is tempting to write

$$\frac{d}{dt} B_t^H = \sum_{i=0}^{\infty} M_H e_i(t) \omega(e_i),$$

which is illegitimate in the traditional sense as the Brownian motion or the fractional Brownian motion is nowhere differentiable. With the chaos expansion form, differentiation and integration with respect to time $t$ can be defined, but they may not always be square integrable. Such type of operation is called integration or differentiation in $(S)^*$ (see [18]).

**Definition 2.1.** Let (S) be a subset of $L^2(\Omega, B(\Omega), \mu)$ consisting of functionals with Wiener-Itô Chaos decomposition such that

$$\sum_{\alpha} \left\{ c_\alpha^2 \alpha! \prod_{j \in N} (2j)^{k\alpha_j} \right\} < \infty$$

for all $k < \infty$ and that $(S)^*$ consists of all expansions, not necessarily belonging to $L^2(\Omega, B(\Omega), \mu)$, such that

$$\sum_{\alpha} \left\{ c_\alpha^2 \alpha! \prod_{j \in N} (2j)^{-q\alpha_j} \right\} < \infty$$

for some $q < \infty$, then, the spaces (S) and $(S)^*$ are called the Hida test function space and the Hida distribution space respectively.

The derivative of the fractional Brownian motion, or the white noise is defined by

$$W^H(t) = \sum_{i=0}^{\infty} M_H e_i(t) \omega(e_i).$$

It can be shown that the sum $W^H(t) \in (S)^*$ (see [14, 31]). The importance of Hida test function space and distribution space is their closedness of Wick multiplication. Wick product is an operator acting on two functionals $F(\omega)$ and $G(\omega)$.

**Definition 2.2.** The Wick's product for two functionals having Wiener-Itô Chaos Decomposition
$$F(\omega) = \sum_\alpha c_\alpha H_\alpha(\omega)$$
and
$$G(\omega) = \sum_\beta b_\beta H_\beta(\omega)$$
is defined as
$$F(\omega) \diamond G(\omega) = \sum_{\alpha,\beta} c_\alpha b_\beta H_{\alpha+\beta}(\omega).$$
Addition of indexes refers to pairwise addition.

The closeness of Wick's product is shown in the following theorem (see Corollary 2.2 and Remark 2.8 of [31]).

**Theorem 2.2.** *Wick multiplication is closed in (S) and $(S)^*$.*

It is reasonable to define the stochastic integration of a functional $Z_t(\omega)$ with respect to the fractional Brownian motion as the integration with respect to time $t$ of the Wick's product between $Z_t(\omega)$ and $W_t^H(\omega)$. Under the Wiener Chaos decomposition framework, if the decomposition exists, the functional $Z_t(\omega) \diamond W_t^H(\omega)$ can be written as
$$\sum_\alpha c_\alpha(t) H_\alpha(\omega).$$
It is natural to think that the integration is
$$\sum_\alpha \left\{ \int_0^t c_\alpha(s) ds \right\} H_\alpha(\omega),$$
by assuming that summation and integration are interchangeable. If the integration with respect to time is a path-wise classical Riemann integral, it is not clear that summation and integration are interchangeable. Such difficulties can be finessed by introducing new definitions for integration with respect to time and integration with respect to the fractional Brownian motion as follows (see Definitions 2.3 and 2.4 respectively).

**Definition 2.3.** (a) Elements in $(S)^*$ as an operator: Let $F(\omega) = \sum_\alpha c_\alpha H_\alpha(\omega) \in (S)^*$ and $f(\omega) = \sum_\alpha b_\alpha H_\alpha(\omega) \in (S)$, then $F$ can be regarded as an operation on $f$
$$\langle F, f \rangle = \sum_\alpha b_\alpha c_\alpha \alpha!.$$

(b) Time integration: If $F_t(\omega) = \sum_\alpha c_\alpha(t) H_\alpha(\omega)$ are elements in $(S)^*$ for all positive real number $t$, and that $\langle F_t, f \rangle$ are integrable with respect to $t$ for all $f \in (S)$, then the integral $\int_R F_t(\omega) dt$ is defined as the unique element in $(S)^*$, $I(\omega)$ such that
$$\langle I(\omega), f \rangle = \int_R \langle F_t(\omega), f \rangle \, dt$$
for all $f \in (S)$.

**Remark 2.2.** It can be shown that the quantity $\langle F, f \rangle$ in part (a) exists and is finite under the condition given in the definition. It can be regarded as the expectation of the product between $F(\omega)$ and $f(\omega)$ when $F(\omega) \in L^2(\Omega, B(\Omega), \mu)$ and $f(\omega) \in (S)$.

**Definition 2.4.** Let $Z(t) = \sum_\alpha c_\alpha(t) H_\alpha(\omega) \in (S)^*$ for any given $t$, then, the Wick's integral of $Z(t)$ is defined as

$$\int_R Z(t) \diamond dB^H(t) = \int_R \{Z(t) \diamond W^H(t)\}\, dt,$$

when $Z(t) \diamond W^H(t)$ is integrable with respect to time $t$ in the sense as in Definition 2.3.

The following theorem asserts that integration and summation are interchangeable, see Lemma 2.5.6 of [18].

**Theorem 2.3.** *Let $Z(t) : R \to (S)^*$, with Wiener Chaos decomposition*

$$\sum_\alpha c_\alpha(t) H_\alpha(\omega)$$

*such that*

$$\sum_\alpha \alpha! \{\int_R c_\alpha(t) dt\}^2 \prod_{j \in N}(2j)^{-q\alpha_j} < \infty$$

*for some $q < \infty$, then $Z(t)$ is time-integrable, also, integration and summation are interchangeable, i.e.*

$$\int_R Z(t) dt = \sum_\alpha \{\int_R c_\alpha(t) dt\} H_\alpha(\omega).$$

## 3. Fractional Itô's lemma

Several approaches of extending classical Itô's lemma to incorporate the fractional Brownian motion were discussed in the literature, e.g., [8, 9, 12]. The settings in these papers are different and various conditions are required to ensure that the stochastic integrals appear in the fractional Itô lemma exist. Bender in [1] provided a simpler version of fractional Itô's lemma based on the white noise setting introduced in the preceding section. Here, we restate Bender's theorem in Theorem 3.1 and give a generalized result in Theorem 3.2.

**Theorem 3.1.** *Consider the stochastic process*

$$Y_t = \int_0^t h(t)\, dB_t^H \equiv \omega[M_H(h(\cdot) 1_{[0,t)}(\cdot))],$$

*where $h(t)$ is a continuous function in $[0, T]$ and $H > \frac{1}{2}$. Let $g(t, y)$ be a two-dimensional function differentiable with respect to $t$ and is twice differentiable with respect to $y$. Also, there exists constants $C_1 \geq 0$ and $\lambda_1 < (2T^H \sup_{s \in [0,T]} h(s))^{-2}$ so that*

$$\max\{|g|, |g_t|, |g_y|, |g_{yy}|\} \leq C_1 e^{\lambda_1 y^2}.$$

Let $k(t) = H(2H-1)h(t)\int_0^T h(s)|s-t|^{2H-2}\,ds$. Then, we have the following fractional Itô's lemma:

$$g(T, Y_T) = g(0,0) + \int_0^T \frac{\partial}{\partial t}g(t, Y_t)\,dt + \int_0^T h(t)\frac{\partial}{\partial y}g(t, Y_t) \diamond dB_t^H$$
$$+ \int_0^T k(t)\frac{\partial^2}{\partial y^2}g(t, Y_t)\,dt.$$

**Remark 3.1.** The integrals above are well defined as the condition that the integrands and integrals both belong to $(L^2)$ is ensured by the given assumption. Since $(L^2) \subset (S)^*$, the integrals are defined as in Section 2.

This theorem gives the differential form of $g(t, Y_t)$ when the underlying stochastic process $Y_t$ is an integrals of a deterministic function $h_t$. The following generalization takes the underlying stochastic process to be $X_t = g(t, Y_t)$ and gives the differential form for $P(t, X_t)$, where $P(t, x)$ is a two-dimensional real-valued function. Before introducing our results, let us illustrate some ideas through an example (see [19]).

Consider the two-dimensional function,

$$g(t, y) = \exp(\mu t - \frac{1}{2}\sigma^2 t^{2H} + \sigma y),$$

where $\mu$ and $\sigma$ are two positive constants and the underlying stochastic process is

$$Y_t = \int_0^t dB_s^H,$$

i.e., $h(t) = 1$. Clearly, for any given value of $T$, the functions $g$, $g_t$, $g_y$ and $g_{yy}$ are all continuous in the closed interval $[0, T]$ and hence, the conditions in the theorem are fulfilled. Applying the lemma, we have

$$k(t) = H(2H-1)\int_0^t |t-s|^{2H-2}\,ds = Ht^{2H-1},$$

and

$$dX_t = (\mu - \sigma^2 Ht^{2H-1})g(t, Y_t)\,dt + \sigma g(t, Y_t) \diamond dB_t^H$$
$$+ k(t)\sigma^2 g(t, X_t)\,dt$$
$$= \mu g(t, Y_t)\,dt + \sigma g(t, Y_t) \diamond dB_t^H$$
$$= \mu X_t\,dt + \sigma X_t \diamond dB_t^H.$$

The next question is whether there exists an Itô's lemma that further expresses $P(t, X_t)$ in terms of integrals involving $\mu X$ and $\sigma X$, but not $Y$ and $g$ explicitly. It is reasonable to expect that

$$dP(t, X_t) = P_t(t, X_t)\,dt + \mu Y_t P_x(t, X_t)\,dt + \sigma Y_t P_x(t, Y_t) \diamond dB_t^H +$$
$$\sigma^2 Ht^{2H-1}X^2 P_{xx}(t, X_t)\,dt.$$

This result can be verified by the following theorem.

**Theorem 3.2.** *Using the same notations and assumptions in Theorem 3.1, further assume that the differential of the stochastic process $X_t = g(t, Y_t)$ can be, according to Theorem 3.1, written as*

$$dX_t = \mu(t, X_t)\,dt + \sigma(t, X_t) \diamond dB_t^H.$$

Also, there exists constants $C_2 \geq 0$ and $\lambda_2 < (2T^H \sup_{s \in [0,T]} h(s))^{-2}$ so that the composite function $P \circ g \equiv P(t, g(t, y))$ satisfies

$$\max\{|(P \circ g)|, |(P \circ g)_t|, |(P \circ g)_y|, |(P \circ g)_{yy}|\} \leq C_2 e^{\lambda_2 y^2}.$$

Let $C(t) = k(t)h^{-2}(t)$. Then the fractional Itô's lemma is given by

$$P(T, X_T) = P(0, X_0) + \int_0^T \frac{\partial}{\partial t} P(t, X_t) \, dt + \int_0^T \mu(t, Y_t) \frac{\partial}{\partial x} P(t, X_t) \, dt$$
$$+ \int_0^T \sigma(t, Y_t) \frac{\partial}{\partial x} P(t, X_t) \diamond dB_t^H + \int_0^T C(t) \sigma^2(t, Y_t) \frac{\partial^2}{\partial x^2} P(t, X_t) \, dt.$$

*Proof.* This theorem can be verified by applying Theorem 3.1 to $f(t, g(t, Y_t))$. □

To relate this result with previous works on fractional stochastic calculus and option pricing theory, such as [12, 26], consider the last term in the right-hand side of the Itô's formula. This second order correction term can be considered as the product of three quantities: $P_{xx}(t, X)$, $\sigma(t, X)$ and $C(t)\sigma(t, X)$. Comparing with the result of [12], the quantity $C(t)\sigma(t, X)$ corresponds to $D_t^\phi X_t$. This quantity is known as the Malliavin derivative of $X_t$. For details of Malliavin calculus, see [20, 24, 27]. Under our assumptions that $X_t$ has the form of $g(t, Y_t)$, this quantity depends only on the current time $t$ and the value of $X$ at time $t$. But this needs not be the situation for a general $X_t$ governed by a fractional stochastic differential equation. In general this quantity may depend on the entire path of $X_t$, not only on the value at one point. This may introduce further complications when working with the differential.

## 4. The fractional CEV model

Here we construct a fractional version of the Constant Elasticity of Variance model by means of the fractional Itô's lemma (Theorem 3.1). As discussed in the introduction, constant elasticity is characterized by the the volatility term $\sigma X_t^{\beta/2}$ in a stochastic diffusion equation. In order to handle long-memory, we replace the Wiener process by a fractional Brownian motion. The fractional diffusion equation is then defined as
$$dX_t = \mu(t, X_t) \, dt + \sigma X_t^{\beta/2} \diamond dB_t^H,$$
where $0 \leq \beta \leq 2$. If $H = \frac{1}{2}$ and $\mu(t, X_t) \equiv \mu X_t$, this is the classical CEV model. In this situation, the Wick integral is equivalent to the Itô's integral (see [18]) and hence, the classical Itô's lemma can be applied to any stochastic process of the form $Y_t = P(t, X_t)$. When long-memory is considered, the Itô's lemma will involve the Malliavin derivative, which is in general path dependent and difficult to handle. When the integration is defined in the white noise sense (Section 2), the integrand is assumed to belong to $(S)^*$ at every time $t$, and the integral is a random variable in $(S)^*$. The elements in $(S)^*$, which are merely formal expansions, may not correspond to real values for each $\omega$ and the term $P(t, X_t)$ may not be well-defined in general. In order to overcome such difficulties, we need to choose a suitable $\mu(t, X_t)$.

Assume that $X_t$ can be written explicitly in terms of time $t$ and a stochastic integral process $Y_t = \int_0^t h(s) dB_s^H$, i.e. $X_t = g(t, Y_t)$. From Theorem 3.1, the differential of $X_t$ can be decomposed into two parts, the drift term and the volatility term. In order to keep the elasticity constant, the function $g(t, y)$ must be chosen as

$$g(t, y) = \left\{ \sigma(1 - \frac{\beta}{2})[h(t)]^{-1} y + f(t) \right\}^{\frac{2}{2-\beta}},$$

where $f(t)$ is an arbitrarily chosen function of $t$. This $g(t, y)$ is in fact the solution to the differential equation

$$h(t) \frac{\partial}{\partial y} g(t, y) = \sigma \{g(t, y)\}^{\beta/2}.$$

The next question is how to choose $f(t)$. The answer is given by the following theorem.

**Theorem 4.1.** *Assume that $h(t)$ is a strictly positive function and $g(t, y)$ is the solution to the differential equation*

(4.1) $$h(t) \frac{\partial}{\partial y} g(t, y) = \sigma(t, g(t, y)).$$

*The general solution to this equation involves an arbitrary function $f(t)$.*

*Let $\eta(t)$ and $\varphi(t)$ be two functions determined by the integral equations,*

$$h(t) = e^{-\int_0^t \eta(s) ds}$$

*and*

$$f(t) = [h(t)]^{-1} \left[ a_0 + \int_0^t h(s) \varphi(s) \, ds \right],$$

*where $a_0$ is a constant so that $g(0, 0) = X_0$, then,*

$$X_t = g(t, Y_t)$$

*yields the volatility $\sigma(t, X_t)$ in the fractional Itô's lemma and the drift term is given by a two dimensional function*

$$\mu(t, x) = \sigma(t, x) \left[ \varphi(t) + C(t) \frac{\partial \sigma}{\partial x} + \int_0^x \frac{\eta(t)}{\sigma(t, x)} dx + \int_0^x \frac{1}{\sigma^2(t, x)} \frac{\partial \sigma}{\partial t} dx \right].$$

*Proof.* Let

$$a(t, x) = \int \frac{dx}{\sigma(t, x)},$$

then $g(t, y)$ given by

$$a(t, g(t, y)) = [h(t)]^{-1} y + f(t)$$

satisfies Equation (4.1). After some manipulations, we have

$$\frac{\partial}{\partial t} g(t, y) + k(t) \frac{\partial^2}{\partial y^2} g(t, y)$$
$$= \sigma(t, g) \left[ \varphi(t) + C(t) \frac{\partial \sigma}{\partial x}(t, g) + \eta(t) a(t, g) - \frac{\partial a}{\partial t} \right],$$

which does not involve $y$ explicitly. The proof is completed by comparing this with fractional Itô's lemma (Theorem 3.1). □

**Remark 4.1.** Note that the function $C(t)$ depends on the choice of $\eta(t)$ and by definition (in Theorem 3.2) must be strictly positive.

In the CEV case, when the drift term has the form of

$$(4.2) \qquad \mu(t,x) = \frac{\eta(t)}{1-\beta/2}x + \sigma\varphi(t)x^{\beta/2} + \frac{\sigma^2\beta}{2}x^{\beta-1}C(t),$$

the solution to the stochastic diffusion equation

$$dX_t = \mu(t, X_t)\,dt + \sigma X_t^{\beta/2} \diamond dB_t^H$$

is well defined and is given by Theorem 4.1. Now choose suitable time dependent functions $\eta(t)$ and $\varphi(t)$. Consider the three quantities of $\mu(t, X_t)$ in (4.2). The last one is strictly positive and cannot be eliminated. When we choose $\eta(t) \equiv (1-\frac{\beta}{2})\mu$ and $\varphi(t) \equiv 0$, the second term vanishes and the first term becomes $\mu X_t$. The drift term becomes

$$\mu(t, X_t) = \mu X_t + \frac{\sigma^2\beta}{2}C(t)X_t^{\beta-1}.$$

The CEV stochastic differential equation is thus

$$X_t = X_0 + \int_0^t \mu X_s\,ds + \left\{ \int_0^t \frac{\sigma^2\beta}{2}C(s)X_s^{\beta-1}\,ds + \int_0^t \sigma X_s^{\beta/2} \diamond dB_s^H \right\}.$$

## 5. Fractional option pricing theory

By using the generalized Itô's lemma (Theorem 3.2), the differential of $P(t, X_t)$ can be decomposed into two terms, the drift one and the volatility one and both terms involve only current time $t$ and $X_t$, i.e., they are path-independent. The foundation of the Black-Scholes' option pricing theory is constructing a self-financing strategy, which makes use of stocks and bonds to hedge an option. The definition for self-financing strategy in continuous-time model depends on how the stochastic integrals are defined. As the fractional stochastic integrals are defined in a different manner, a new definition for self-financing strategy is required. One approach is adopted by [14, 19], which defines self-financing strategy under the geometric fractional Brownian motion. Here, this approach is extended to a more general situation.

Let $X_t$ be the stock price process and $\Pi_t$ be the bond value and they are governed by

$$X_t = \int_0^t \mu(s, X)\,ds + \int_0^t \sigma(s, X) \diamond dB_s^H,$$
$$(5.1) \qquad \Pi_t = \int_0^t r\Pi\,ds.$$

**Definition 5.1.** (see Section 5 of [14]) A trading strategy consists of a quantity $(u_t, v_t)$ of bonds and stocks is called self-financing if the infinitesimal change in the portfolio value at time $t$ is given by

$$\begin{aligned} dZ_t &= d(u_t\Pi_t + v_tX_t) \\ &= r\Pi_t u_t\,dt + \mu(t, X_t)v_t\,dt + [\sigma(t, X_t)v_t] \diamond dB_t^H + d\Delta, \end{aligned}$$

where $d\Delta$ is an infinitesimal dividend payment term.

Below, using the above definition, a fractional Black-Scholes equation is derived.

**Theorem 5.1.** *Suppose that the market consists of two securities, a risk-free bond and a stock. Here, the stock provides dividend continuously with rate $\delta$. Assume that the stock price process $X_t = g(t, Y_t)$ is defined as in Section 3 which also satisfies the equations (5.1). Then the price of a derivative on the stock price with a bounded payoff $f(X(T))$ is given by $P(t, X)$, where $P(t, X)$ solves the PDE:*

$$(5.2) \qquad \frac{\partial P}{\partial t} + \sigma^2(t, X)C(t)\frac{\partial^2 P}{\partial X^2} + (r - \delta)X\frac{\partial P}{\partial X} - rP = 0,$$

*with boundary condition $P(T, X) = f(X)$ given that $P \circ g$ satisfies the conditions in Theorem 3.2.*

*Proof.* A proof parallel to the fractional Black-Scholes equation in [26] is given here. Consider a solution $P(t, X)$ to equation (5.2). Applying fractional Itô's Lemma Theorem 3.2,

$$dP(t, X_t) = [\frac{\partial P}{\partial t} + \frac{\partial P}{\partial X}\mu(t, X) + \frac{\partial^2 P}{\partial X^2}\sigma^2(t, X)C(t)]\, dt + \frac{\partial P}{\partial X}\sigma(t, X) \diamond dB_t^H.$$

Form a trading strategy by dynamically adjusting a portfolio consisting a varying quantity $v(t)$ of stocks and $u(t)$ of bonds. By choosing

$$v(t) = \frac{\partial P}{\partial X}$$
$$(5.3) \qquad u(t) = \frac{1}{\Pi_t}(P - X\frac{\partial P}{\partial X}),$$

then the portfolio value at time $t$ is $P_t$ and

$$ru(t)\Pi_t\, dt + v(t)\mu(t, X_t)\, dt + [v(t)\sigma(t, X_t)] \diamond dB_t^H + \delta v(t)X_t\, dt$$
$$= [rP - rX\frac{\partial P}{\partial t}]\, dt + \mu\frac{\partial P}{\partial X}\, dt + (\sigma\frac{\partial P}{\partial X}) \diamond dB_t^H + \frac{\partial P}{\partial X}\delta X\, dt$$
$$= [\frac{\partial P}{\partial t} - \frac{\partial P}{\partial X}\delta X\, dt + \frac{\partial^2 P}{\partial X^2}\sigma^2(t, X)C(t)]\, dt$$
$$\qquad + \mu\frac{\partial P}{\partial X}\, dt + (\sigma\frac{\partial P}{\partial X}) \diamond dB_t^H + \frac{\partial P}{\partial X}\delta X\, dt$$
$$= [\frac{\partial P}{\partial t} + \frac{\partial P}{\partial X}\mu(t, X) + \frac{\partial^2 P}{\partial X^2}\sigma^2(t, X)C(t)]\, dt + \frac{\partial P}{\partial X}\sigma(t, X) \diamond dB_t^H$$
$$= dP(t, X_t)$$
$$= d(u_t\Pi_t + v_tX_t).$$

By Definition 5.1, $(u(t), v(t))$ is a self-financing strategy. It can be shown that such strategy hedges the derivative. The portfolio value at time $t$ is given by $u(t)\Pi_t + v(t)X_t$ and it is equal to $P(t, X_t)$. At time of maturity, the portfolio value is just $P(T, X_T)$. By assumption, the function $P(t, X)$ satisfies the boundary condition, so $P(T, X_T) = f(X_T)$. Therefore $(u(t), v(t))$ hedges the derivative and $P(t, X)$ is the option price. □

## 6. Pricing an European call option under CEV models

Putting $\sigma^2(t, X) \equiv \sigma^2 X^\beta$, the Black-Scholes PDE (Theorem 5.1) of the CEV model is now given by

$$\frac{\partial P}{\partial t} + \sigma^2 C(t)X^\beta\frac{\partial^2 P}{\partial X^2} + (r - \delta)X\frac{\partial P}{\partial X} - rP = 0.$$

Putting $Y = X^{2-\beta}$ (see [5]) and $P(t, X) = e^{rt}Q(t, Y)$, this equation becomes

$$\frac{\partial Q}{\partial t} + [bY + cC(t)]\frac{\partial Q}{\partial Y} + aC(t)Y\frac{\partial^2 Q}{\partial Y^2} = 0,$$

where $a = \sigma^2(2-\beta)^2$, $b = (r-\delta)(2-\beta)$ and $c = \sigma^2(2-\beta)(1-\beta)$. The boundary condition is

$$Q(T, Y) = e^{-rT}\max(Y^{\frac{1}{2-\beta}}, 0).$$

The approach of Cox and Ross [6] that made use of Feller's result ([15, 16]) can be adopted. First the solution for this equation at $(t_0, Y_{t_0})$ is the expectation of $Q(T, Y_T)$ under the SDE

(6.1) $$dY = [bY + cC(t)]\,dt + \sqrt{2aC(t)Y}\,dB_t,$$

with $Y(t_0) = Y_{t_0}$ (see [25]). To solve this SDE, we follow Feller's arguments. First, a useful result of Kolmogorov is stated below (see Equation (167) of [21]).

**Theorem 6.1.** *The probability density function of a diffusion process $X_t$ driven by standard Brownian motion*

$$dX_t = \mu(t, X_t)\,dt + \sigma(t, X_t)dB_t$$

*is given by the PDE*

$$u_t = [\frac{1}{2}\sigma^2(t, X)u(t, X)]_{XX} - [\mu(t, X)u(t, X)]_X.$$

In our case, because of (6.1), the Kolmogorov's equation becomes

$$u_t = [aC(t)Yu(t, Y)]_{YY} - [(bY + cC(t))u(t, Y)]_Y.$$

The European call option pricing formula can be obtained by solving the above PDE.

**Theorem 6.2.** *Under the fractional CEV model introduced in Section 4, the price of an European call option with strike price $K$, mature at $T$ at current time $t_0$ is given by*

$$P(t_0, X_0) = e^{-\delta(T-t_0)}X_0 \sum_{r=0}^{\infty} \frac{1}{r!} e^{-\frac{x}{a\gamma_T}} (\frac{x}{a\gamma_T})^r G(r+1+\frac{1}{2-\beta}, \frac{K^{2-\beta}}{ae^{bT}\gamma_T})$$

$$- Ke^{-r(T-t_0)} \sum_{r=0}^{\infty} \frac{1}{\Gamma(r+1-\frac{1}{2-\beta})} e^{-\frac{x}{a\gamma_T}} (\frac{x}{a\gamma_T})^{r+\frac{1}{2-\beta}} G(r+1, \frac{K^{2-\beta}}{ae^{bT}\gamma_T}),$$

*where* $x = e^{-bt_0}Y(t_0)$,

$$\gamma_t \equiv \int_{t_0}^{t} e^{-b\tau}C(\tau)\,d\tau$$

*and*

$$G(\alpha, \nu) \equiv \frac{1}{\Gamma(\alpha)} \int_{\nu}^{\infty} e^{-\tau}\tau^{\alpha-1}\,d\tau.$$

*Proof.* Assume that the Laplace Transform of $u(t, Y)$ with respect to $Y$ exists and equals to $\omega(t, s)$. Since the value of $Y$ at time $t_0$ is given, $Y_{t_0}$ is deterministic and thus $u(t_0, Y) = \delta(Y - Y_{t_0})$, the Dirac function. Also, $\omega(t_0, s) = e^{-sY_0} = \pi(s)$ and the equation becomes the boundary value problem

(6.2) $$\omega_t + s(aC(t)s - b)\omega_s = -cC(t)s\omega + \psi(t)$$
(6.3) $$\omega(t_0, s) = e^{-sY_0}.$$

In equation (6.2), $\psi(\cdot)$ is called the flux of $u$ at the origin (see [16]), which is to be determined later. Now, we find the characteristic curve of the first order PDE (6.2). The characteristic curve is given by

$$\frac{ds}{dt} = aC(t)s^2 - bs.$$

The solution to this equation is

$$s = e^{-bt}[C_1 - a\gamma_t]^{-1}.$$

On this curve, $\omega(t, s(t))$ satisfies

$$\frac{d}{dt}\omega(t, s(t)) = \psi(t) - cC(t)s(t)\omega(t, s(t)).$$

Solving, we have

$$\omega(t, s(t)) = [C_1 - a\gamma_t]^{c/a}[C_2 + \int_{t_0}^{t} \psi(\tau)|C_1 - a\gamma_\tau|^{-c/a} d\tau].$$

For any given point $(t_1, s_1)$, the characteristic curve with

$$C_1 = a\gamma_{t_1} + \frac{1}{s_1 e^{bt_1}} = C(t_1, s_1)$$

passes through $(t_1, s_1)$. Also, this curve passes through the point $(t_0, C_1^{-1}e^{-bt_0})$. This yields the value of $C_2$,

$$C_2 = [C_1(t_1, s_1)]^{-c/a}\omega(t_0, e^{-bt_0}C_1^{-1}(t_1, s_1)).$$

The Laplace transform of $u(t, Y)$ at point $(t_1, s_1)$ is thus given by

(6.4) $$\begin{aligned}\omega(t_1, s_1) &= (s_1 e^{bt_1})^{-c/a}[(C_1(t_1, s_1))^{-c/a}\omega(t_0, e^{-bt_0}C_1^{-1}(t_1, s_1)) \\ &\quad + \int_{t_0}^{t_1} \psi(\tau)|a(\gamma_{t_1} - \gamma_\tau) + \frac{1}{s_1 e^{b\tau}}|^{-c/a} d\tau] \\ &= [s_1 a e^{bt_1}\gamma_{t_1} + 1]^{-c/a}\pi(\frac{s_1 e^{b(t_1-t_0)}}{s_1 a e^{bt_1}\gamma_{t_1} + 1}) \\ &\quad + \int_{t_0}^{t_1} [s_1 a e^{bt_1}(\gamma_{t_1} - \gamma_\tau) + 1]^{-c/a}\psi(\tau) d\tau.\end{aligned}$$

Following the argument of [16], when $u(t, 0)$ is finite and $c \leq 0$ or $0 < c < a$,

$$\lim_{s \to}(sae^{bt}\gamma_t + 1)\omega(t, s) = 0,$$

then $\psi(t)$ is given by the integral equation

$$\pi(\frac{e^{-bt_0}}{a\gamma_t}) + \int_{t_0}^{t} \psi(\tau)(\frac{\gamma_t}{\gamma_t - \gamma_\tau})^{c/a} d\tau = 0.$$

To solve this equation, applying the substitutions $z = \frac{1}{\gamma_t}$ and $\zeta = \frac{1}{\gamma_\tau}$,

$$\int_z^\infty g(\zeta)(\zeta - z)^{-c/a}\, d\zeta = \pi\left(\frac{e^{-t_0}}{a\gamma_t}\right) = \pi\left(\frac{ze^{-t_0}}{a}\right)$$
$$= \exp\left(-\frac{xz}{a}\right),$$
(6.5)

where $g(\zeta) = -\psi(\tau)\zeta^{c/a}\frac{d\tau}{d\zeta}$.

The solution of (6.5) is

$$g(\zeta) = \frac{1}{\Gamma(1 - \frac{c}{a})}\left(\frac{x}{a}\right)^{\frac{-c+a}{a}}\exp\left(-\frac{x}{a\gamma_\tau}\right),$$

$$\psi(\tau) = g(\zeta)\zeta^{-c/a}\frac{d\zeta}{d\tau}$$
$$= \frac{-1}{\Gamma(1 - \frac{c}{a})}\left(\frac{1}{\gamma_\tau}\right)\left(\frac{x}{a\gamma_\tau}\right)^{\frac{-c+a}{a}}\exp\left(-\frac{x}{a\gamma_\tau}\right)\frac{d\gamma_\tau}{d\tau}.$$

Substituting this result into (6.4), after simplification, we get

$$\omega(t, s) = \exp\left(\frac{-sxe^{bt}}{sae^{bt}\gamma_t + 1}\right)\left(\frac{1}{sae^{bt}\gamma_t + 1}\right)^{c/a}\Gamma\left(1 - \frac{c}{a}; \frac{x}{a\gamma_t(sae^{bt}\gamma_t + 1)}\right).$$

The next step is to perform an inverse Laplace transform with respect to $s$. To this end, let

$$A = \frac{x}{a\gamma_t},$$
$$z = sae^{bt}\gamma_t + 1.$$

One can verify that equation (6.5) in [16] is still valid after these substitutions. The quantity $\omega(t, s)$ can now be rewritten as

$$\frac{1}{\Gamma(1 - \frac{c}{a})}e^{-A}A^{1-\frac{c}{a}}\int_0^1 (1 - \tau)^{-c/a}e^{\frac{A\tau}{z}}z^{-1}\, d\tau.$$

Using the fact that Laplace Transform of $I_0(z(A\tau Y)^{1/2})$ is $e^{\frac{A\tau}{z}}z^{-1}$, we have

$$u(t, Y) = \left(\frac{1}{ae^{bt}\gamma_t}\right)\left(\frac{xe^{bt}}{Y}\right)^{\frac{-c+a}{2a}}\exp\left\{-\frac{(Y + xe^{bt})}{ae^{bt}\gamma_t}\right\}I_{1-\frac{c}{a}}\left[\frac{2}{a\gamma_t}(e^{-bt}xY)^{1/2}\right],$$

where $I_\lambda(\cdot)$ is the first type Bessel function with order $\lambda$, which is defined as

$$I_\lambda(\cdot) = \sum_{k=0}^\infty \frac{(\cdot/2)^{2k+\lambda}}{k!\Gamma(k + 1 + \lambda)}.$$

This density function is then used to find the solution of $P$ at $(t_0, X_0)$ by means of the identity,

$$P(t_0, X_0) = e^{rt_0}Q(t_0, X_0) = e^{-r(T-t_0)}E[\max(Y_T^{\frac{1}{2-\beta}} - K, 0)].$$

After direct calculations, we have

$$P(t_0, X_0) = e^{-r(T-t_0)}\int_{K^{2-\beta}}^\infty (y^{\frac{1}{2-\beta}} - K)u(T, y)\, dy$$
$$= e^{-r(T-t_0)}\sum_{r=0}^\infty e^{(r-\delta)T}x^{\frac{1}{2-\beta}}\frac{1}{r!}e^{-\frac{x}{a\gamma_T}}\left(\frac{x}{a\gamma_T}\right)^r G\left(r + 1 + \frac{1}{2-\beta}, \frac{K^{2-\beta}}{ae^{bT}\gamma_T}\right)$$
$$- Ke^{-r(T-t_0)}\sum_{r=0}^\infty \frac{1}{\Gamma(r + 1 - \frac{1}{2-\beta})}e^{-\frac{x}{a\gamma_T}}\left(\frac{x}{a\gamma_T}\right)^{r+\frac{1}{2-\beta}} G\left(r + 1, \frac{K^{2-\beta}}{ae^{bT}\gamma_T}\right).$$

Putting $x = e^{-bt_0} X_0^{2-\beta}$,

$$x^{\frac{1}{2-\beta}} = e^{-\frac{bt_0}{2-\beta}} X_0 = e^{-(r-\delta)t_0} X_0.$$

So the option pricing formula is

$$P(t_0, X_0) = e^{-\delta(T-t_0)} X_0 \sum_{r=0}^{\infty} \frac{1}{r!} e^{-\frac{x}{a\gamma_T}} \left(\frac{x}{a\gamma_T}\right)^r G\left(r+1+\frac{1}{2-\beta}, \frac{K^{2-\beta}}{ae^{bT}\gamma_T}\right)$$
$$- Ke^{-r(T-t_0)} \sum_{r=0}^{\infty} \frac{1}{\Gamma(r+1-\frac{1}{2-\beta})} e^{-\frac{x}{a\gamma_T}} \left(\frac{x}{a\gamma_T}\right)^{r+\frac{1}{2-\beta}} G\left(r+1, \frac{K^{2-\beta}}{ae^{bT}\gamma_T}\right).$$

$\square$

This formula is similar to the classical one, which is obtained by replacing the term $\gamma_T$ by the term $\frac{1}{2b}(e^{-bt_0} - e^{-bT})$. As these two terms do not depend on the strike price, the implied volatility pattern is the same as the classical CEV model. Consequently, the fractional CEV model also accounts for the volatility skewness observed in practice.